\documentclass[12pt]{article}
\usepackage{amsfonts}
\usepackage{amsmath}
\usepackage{amssymb}
\usepackage{amsthm}

\newtheorem{prop}{Proposition}
\newtheorem{lemma}{Lemma}
\newtheorem{theorem}{Theorem}
\newtheorem{rem}{Remark}
\newtheorem{exmp}{Example}
\newtheorem*{cor}{Corollary}

\renewcommand{\thepage}{\roman{page}}

\title{Transformations of Grassmann spaces I\\
(preliminary version)}
\author{Mark  Pankov}
\date{}

\begin{document}
\maketitle
\tableofcontents
\newpage
\section*{Preface}

This is a preliminary version of the first part of the book
"Transformations of Grassmann spaces".

Let $\Pi$ be an $n$-dimensional projective space and $n\ge 3$.
If $0<k<n-1$ then the Grassmann space ${\mathcal G}_{k}(\Pi)$
has the natural structure of a partial linear space such that
two points are collinear if and only if they are adjacent subspaces.
For this case the classical Chow's theorem \cite{Chow} says that
any bijective transformation of ${\mathcal G}_{k}(\Pi)$
preserving the adjacency relation in both directions is induced by
a collineation of $\Pi$ to itself or to the dual space $\Pi^{*}$,
the second possibilities can be realized only for the case when $n=2k+1$.
In Section 2 we consider some results closely related with Chow's theorem
\cite{Brauner2}, \cite{Havlicek3}, \cite{Havlicek5}, \cite{Huang1},
\cite{Kreuzer2}.

Let $B$ be a base for $\Pi$.
The set of all $k$-dimensional subspaces spanned by points of $B$
is called the {\it base subset} of ${\mathcal G}_{k}(\Pi)$ associated with $B$.
We prove that mappings of ${\mathcal G}_{k}(\Pi)$ to ${\mathcal G}_{k}(\Pi')$
($\Pi'$ is another $n$-dimensional projective space)
sending base subsets to base subsets are induced by strong embeddings of
$\Pi$ to $\Pi'$ or $\Pi'^{*}$ (Theorem 2.2);
some partial cases of this statement can be found in
\cite{Pankov1}, \cite{Pankov2} and \cite{Pankov3};
a weak version of our result holds for
Grassmann spaces of linear spaces satisfying the exchange axiom
(Theorems 2.3 and 2.4).
The example constructed in \cite{Kreuzer2} shows that
for adjacency preserving mappings of ${\mathcal G}_{k}(\Pi)$ to ${\mathcal G}_{k}(\Pi')$
the same fails.

Base subsets of ${\mathcal G}_{k}(\Pi)$ can be considered as
"levels" of apartments of the building associated with $\Pi$
(this building consists of all maximal flags of $\Pi$,
apartments
are the sets of all maximal flags spanned by points of bases for $\Pi$,
see \cite{Buekenhout}, \cite{Taylor} or \cite{Tits}).
It follows from results given in \cite{AbramenkoVanMaldeghem}
that transformations of a spherical building preserving the class of apartments
are automorphisms of this building.
By Theorem 2.2, something like this holds for Grassmann spaces.

The second part of the book will be devoted to
analogues  of Chow's Theorem and Theorem 2.2 for
polar spaces \cite{Chow}, \cite{Huang2}, \cite{Huang3} \cite{Pankov4}.

\newpage
\renewcommand{\thepage}{\arabic{page}}
\setcounter{page}{1}

\section{Linear spaces}
We start with the definition and a few simple examples of linear spaces.
After that we will study only linear spaces satisfying the exchange axiom;
an important partial case of these spaces are well-known affine and projective spaces
(affine spaces are not considered here).

General properties of morphisms of linear spaces
(semicollineations and embeddings) will be given in Subsection 1.4.
Morphisms of the projective spaces associated with vector spaces
are induced by semilinear mappings (Subsection 1.5).

\subsection{Main definitions}
Let $P$ be a set of points and ${\mathcal L}$ be a family of proper subsets of $P$;
elements of this family will be called {\it lines}.
We say that points $p,q,\dots$ are {\it collinear}
if there exists a line containing them;
otherwise, these points are said to be {\it non-collinear}.
We will suppose that the pair ${\Pi}=(P,{\mathcal L})$ is a {\it linear space};
this means that the following axioms hold true:
\begin{enumerate}
\item[(L1)]
Each line contains at least two points.
\item[(L2)]
For any two distinct points $p_{1}$ and $p_{2}$ there is unique line containing them;
this line will be denoted by $p_{1}p_{2}$.
\end{enumerate}
A set $S\subset P$ is said to be a {\it subspace} of the linear space ${\Pi}$
if for any two distinct points $p_{1}$ and $p_{2}$ belonging to $S$ the line
$p_{1}p_{2}$ is contained in $S$.
By this definition, the empty set and one-point sets are subspaces
(since these sets do not contain two distinct points).
A direct verification shows that
the intersection of any collection of subspaces is a subspace.

Let $X$ be a subset of $P$.
The minimal subspace containing $X$ (the intersection of all subspaces containing $X$)
is called {\it spanned} by $X$ and denoted by $\overline{X}$.
For the case when $\overline{X}$ coincides with $P$
we say that our linear space is spanned by the set $X$.

If our set $P$ contains at least two points then
we put $[X]_{1}$ for the union of all lines $p\,q$ such that
$p$ and $q$ are distinct points of the set $X$
($[X]_{1}:=X$ if $X$ is empty or a one-point set)
and define
$$[X]_{i}:=[[X]_{i-1}]_{1}$$
for each natural number $i \ge 2$;
we will also assume that $[X]_{0}$ coincides with $X$.
It is trivial that $[X]_{i}$ is contained in $[X]_{j}$ if $i<j$.

\begin{lemma}
$\overline{X}=\bigcup^{\infty}_{i=0}\,[X]_{i}$.
\end{lemma}

\begin{proof}
Since each $[X]_{i}$ is contained in the subspace $\overline{X}$,
we need to prove that the union of all $[X]_{i}$ is a subspace.
Any two distinct points $p\in [X]_{i}$ and $q\in [X]_{j}$
belong to the set $[X]_{k}$ where $k=\max\{i,j\}$;
hence the line $p\,q$ is contained in $[X]_{k+1}$.
\end{proof}

A set $X\subset P$ is  said to be {\it independent} if
the subspace $\overline{X}$ is not spanned by a proper subset of $X$
(we will say that points $p,q,\dots$ are independent
if the set formed by them is independent).
A subset of an independent set is independent.

Let $S$ be a subspace of $\Pi$.
An independent subset of $S$ is called a {\it base} for $S$
if $S$ is spanned by it; for the case when $S$ coincides with $P$
we say that this set is a base for our linear space.

In this book we will consider only finite-dimensional linear
spaces (linear spaces spanned by finite sets of points) and say
that a linear space (a subspace of a linear space) is
$n$-{\it dimensional} if $n+1$ is the smallest number of points spanning it.

It follows from the definition of the dimension that
any base for an $n$-dimensional linear space
contains at least $n+1$ points and there exists a base consisting of exactly $n+1$ points.
It is natural to ask: have any two bases the same cardinal number?
We will discuss this question in Subsection 1.3.

The empty set is unique $(-1)$-dimensional subspace.
Points and lines are $0$-dimensional and $1$-dimensional subspaces
(respectively).
A linear space is spanned by at least three non-collinear points
and its dimension is not less than $2$.
In what follows $2$-dimensional linear spaces and
$2$-dimensional subspaces of linear spaces
will be called {\it planes}.

Now we consider two examples.
The first is trivial, the construction given in Example 1.2 will be often
exploited in follows.

\begin{exmp}{\rm
Let $P$ be a finite set of points such that $|P|\ge 3$ and
${\mathcal L}$ be the family of all subsets of $P$
containing exactly two points.
Then $(P,{\mathcal L})$ is a linear space.
Any subset of $P$ is a subspace and the dimension of this linear space is equal to $|P|-1$.
}\end{exmp}

\begin{exmp}{\rm
Let $\Pi =(P, {\mathcal L})$ be a linear space and $S$
be a subspace of $\Pi$ the dimension of which is not less than $2$.
Denote by ${\mathcal L}_{S}$ the family of all lines contained in $S$.
It is trivial that $\Pi_{S}=(S,{\mathcal L}_{S})$ is a linear space.
Now suppose that $X$ is a subset of $P$
containing at least three non-collinear points and
put ${\mathcal L}_{X}$ for the family of all sets $L'\subset X$
satisfying the following conditions:
$|L'|\ge 2$ and there exists a line $L\in {\mathcal L}$
such that
$$L'=L\cap X.$$
Then $\Pi_{X}=(X,{\mathcal L}_{X})$ is a linear space.
}\end{exmp}

Let $\Pi=(P,{\mathcal L})$ and $\Pi'=(P',{\mathcal L}')$ be linear spaces.
A bijection $f:P\to P'$ is called a {\it collineation} of $\Pi$ to $\Pi'$ if
$$f({\mathcal L})={\mathcal L}';$$
in other words,
the bijection $f$ is {\it collinearity} and {\it non-collinearity} preserving
($f$ transfers triples of collinear and non-collinear points
to collinear and non-collinear points, respectively).
The linear spaces $\Pi$ and $\Pi'$ are said to be
{\it isomorphic} if collineations of $\Pi$ to $\Pi'$ exist.

\begin{exmp}{\rm
If all lines of the linear $\Pi$ and $\Pi'$ have only two points
(Example 1.1) then any bijection of $P$ to $P'$ is a collineation.
}\end{exmp}

\begin{lemma}
If $f:P\to P'$ is a collineation then the following statements are fulfilled:
\begin{enumerate}
\item[{\rm(1)}]
$f^{-1}$ is a collineation;
\item[{\rm(2)}]
a set $S\subset P$ is a subspace of $\Pi$ if and only if $f(S)$ is a subspace of $\Pi'$;
\item[{\rm(3)}]
for any set $X\subset P$ we have
$$f(\overline{X})=\overline{f(X)};$$
thus $X$ is independent or a base for $\Pi$
if and only if $f(X)$ is independent or a base for $\Pi'$, respectively;
\item[{\rm(4)}]
a subspace of $\Pi$ and its $f$-image have the same dimension;
in particular, our linear spaces have the same dimension.
\end{enumerate}
\end{lemma}

\begin{proof}
The statements (1) and (2) are trivial,
(3) follows from (2) and the statement
(4) is a direct consequence of (3).
\end{proof}

\subsection{Projective spaces}

First of all we remark that
{\it if any two lines of a linear space have a non-empty intersection
then this space is a plane};
thus if the dimension of a linear space
is not less than $3$ then non-intersecting lines exist.

\begin{proof}
Let $(P,{\mathcal L})$ be a linear space and $p_{1},p_{2},p_{3}\in P$ be non-collinear points.
By the hypothesis, for any point $p\in P-\{p_{3}\}$
the lines $pp_{3}$ and $p_{1}p_{2}$ have a non-empty intersection.
This means that our linear space is spanned by $p_{1},p_{2},p_{3}$.
\end{proof}

Recall that a linear space is called a {\it projective plane} if
the following axioms hold:
\begin{enumerate}
\item[(P1)] any two lines have a non-empty intersection,
\item[(P2)] each line contains at least three points.
\end{enumerate}
Then a {\it projective space} can be defined as
a linear space where each plane is projective.
The classical axioms for projective spaces can be found, for example, in \cite{VeblenYoung}.

\begin{rem}{\rm
A. Kreuzer \cite{Kreuzer3} has discussed the following question:
how many planes of a linear space must be projective to
ensure that this space is projective?
}\end{rem}

\begin{exmp}{\rm
Let $V$ be a left vector space over a division ring.
Put $P(V)$ for the set of all $1$-dimensional subspaces of $V$.
A subset of $P(V)$ will be called a {\it line} if
it consists of all $1$-dimensional subspaces contained in some
$2$-dimensional subspace of $V$.
If $\dim V\ge 3$ then the following statements are fulfilled
(see, for example, \cite{Buekenhout}):
\begin{enumerate}
\item[(1)]
The set $P(V)$ together with the family of lines defined above
is a projective space which dimension is equal to $\dim V -1$;
we will denote this space by $\Pi (V)$.
\item[(2)]
A subset of $P(V)$ is a $k$-dimensional subspace of $\Pi (V)$
if and only if it consists of all $1$-dimensional subspaces
contained in some $(k+1)$-dimensional subspace of $V$.
\item[(3)]
Elements of the set $P(V)$ are independent points of
the projective space $\Pi(V)$
or form a base for it if and only if non-zero vectors lying on these subspaces
are linearly independent or form a base for $V$.
\end{enumerate}
Projective spaces over fields are said to be {\it Pappian};
there is a geometrical characterization of the commutativity
known as the {\it Pappian axiom} (see \cite{Artin} or \cite{Baer}).
By J. H. M. Wedderburn's Theorem \cite{Wedderburn} (see also \cite{Artin}),
any finite division ring is a field.
Projective spaces over finite fields are studied in \cite{Buekenhout},
\cite{Hirschfeld1}, \cite{Hirschfeld2},
\cite{HirschfeldThas}, \cite{Thas}.
}\end{exmp}

It is well-known that (see, for example, \cite{Baer})
{\it if the dimension of a projective space is not less than three
then it is isomorphic to the projective space
associated with certain left vector space over a division ring.}
For a projective plane this statement holds
only for the case when it satisfies the Desarguesian axiom;
note that non-Desarguesian projective planes exist.

\subsection{Exchange axiom}
Let $\Pi=(P,{\mathcal L})$ be a linear space.
We say that $\Pi$ satisfies the {\it exchange axiom}
if for any set $X\subset P$ and for any two points
$p_{1}$ and $p_{2}$ belonging to $P-\overline{X}$
$$p_{2}\in \overline{X \cup \{p_{1}\}}\Longrightarrow p_{1}\in \overline{X\cup \{p_{2}\}}.$$

\begin{exmp}{\rm
If $\Pi$ is  projective  then for any set $X\subset P$
and any point $p$ belonging to $P-\overline{X}$ we have
$$\overline{X\cup \{p\}}=[\overline{X}\cup\{p\}]_{1}$$
(it is a simple consequence of the axiom (P1)).
This implies the fulfillment of the exchange axiom.
}\end{exmp}

\begin{exmp}{\rm
Suppose that $\Pi$ is a projective space and take a point $p\in P$.
By the axiom (P2), the set
$$X:=P-\{p\}$$
contains at least three non-collinear points and
the linear space $\Pi_{X}$ is well-defined (Example 1.2).
A set $S\subset X$ is a subspace of $\Pi_{X}$
if and only if $S$ or $S\cup \{p\}$ is a subspace of $\Pi$.
This means that $\Pi_{X}$ satisfies the exchange axiom
(since the exchange axiom holds for $\Pi$).
Let $L_{1},L_{2}\in {\mathcal L}$ be distinct lines passing through the point $p$.
Then $L_{1}-\{p\}$ and $L_{2}-\{p\}$ are non-intersecting elements of ${\mathcal L}_{X}$;
these lines span a plane of $\Pi_{X}$ and the axiom (P1) does not hold for this space.
Similarly, for any proper subspace $S$ of $\Pi$ the linear space $\Pi_{P-S}$
has the same properties.
}\end{exmp}

\begin{theorem}
If $\Pi$ satisfies the exchange axiom
then any two bases for $\Pi$ have the same cardinal number
and any independent subset of
$P$ is contained in certain base for $\Pi$.
\end{theorem}

To prove Theorem 1.1 we will use the following lemma.

\begin{lemma}
If $\Pi$ satisfies the exchange axiom and
$X=\{p_{1},\dots, p_{m}\}$ is an independent subset of $P$ then
$$\overline{X-\{p_{i_{1}},\dots, p_{i_{p}}\}}\cap
\overline{X-\{p_{j_{1}},\dots, p_{j_{q}}\}}$$
coincides with
$$\overline{X-(\{p_{i_{1}},\dots, p_{i_{p}}\}\cup
\{p_{j_{1}},\dots, p_{j_{q}}\})}.$$
\end{lemma}

\begin{proof}
First of all we show that
$$
\overline{X-\{p_{i}\}}\cap \overline{X-\{p_{j}\}}
=\overline{X-\{p_{i},p_{j}\}}.
$$
Clearly,
$$\overline{X-\{p_{i},p_{j}\}}\subset
\overline{X-\{p_{i}\}}\cap \overline{X-\{p_{j}\}}.$$
If the inverse inclusion fails then
we take a point $p$ belonging to
$$(\overline{X-\{p_{i}\}}\cap \overline{X-\{p_{j}\}})-\overline{X-\{p_{i},p_{j}\}}.$$
By the exchange axiom,
$$\overline{X-\{p_{i}\}}=\overline{(X-\{p_{i},p_{j}\})\cup \{p\}}=
\overline{X-\{p_{j}\}}.$$
This contradicts to the independence of the set $X$
and we get the required.

The equality
$$\overline{X-\{p_{i_{1}}\}}\cap\dots\cap \overline{X-\{p_{i_{k}}\}}=
\overline{X-\{p_{i_{1}},\dots, p_{i_{k}}\}}$$
can be proved by induction;
our statement follows from it.
\end{proof}

\begin{proof}[Proof of Theorem 1.1]
Suppose that $q_{1},\dots,q_{m}$ are independent points of $P$.
Denote by $n$ the dimension of $\Pi$ and consider the base
$B=\{p_{1},\dots, p_{n+1}\}$ for $\Pi$.
By Lemma 1.3,
$$\overline{B-\{p_{1}\}}\cap\dots\cap \overline{B-\{p_{n+1}\}}
=\emptyset.$$
Hence there exists a number $i$ such that
the subspace $\overline{B-\{p_{i}\}}$ does not contain the point $q_{1}$.
Clearly, we can assume that $i=1$.
Then the exchange axiom guarantees that $\{q_{1},p_{2},\dots, p_{n+1}\}$
is a base for our linear space.

Now suppose that for certain number $k\ge 1$ the set
$$B'=\{q_{1},\dots,q_{k},p_{k+1},\dots,p_{n+1}\}$$
is a base for $\Pi$.
Since
$$\overline{B'-\{p_{k+1}\}}\cap\dots\cap \overline{B'-\{p_{n+1}\}}
=\overline{\{q_{1},\dots,q_{k}\}}$$
(Lemma 1.3),
the point $q_{k+1}$ does not belong to some subspace $\overline{B'-\{p_{j}\}}$
such that $j\ge k+1$;
we will assume that $j=k+1$.
By the exchange axiom, the points
$$q_{1},\dots,q_{k+1},p_{k+2},\dots,p_{n+1}$$
form a base for $\Pi$.

Therefore,
our collection of independent points can be extended to a base for
$\Pi$ if $m\le n+1$.
If $m>n+1$ then the arguments given above show that $q_{1},\dots,q_{n+1}$ form a base
for $\Pi$ and the points $q_{1},\dots,q_{m}$ are not independent.
Thus $m\le n+1$.
If $\{q_{1},\dots,q_{m}\}$ is a base for our linear space then $m$ is not less than $n+1$
(by the definition of the dimension),
hence $m=n+1$.
\end{proof}

\subsection{Morphisms of linear spaces: Semicollineations and embeddings}
Let $\Pi=(P,{\mathcal L})$ and $\Pi'=(P',{\mathcal L}')$ be linear spaces.
A bijection $f:P\to P'$ is called a {\it semicollineation} of $\Pi$ to $\Pi'$
if it is collinearity preserving;
in other words, for any line $L\in {\mathcal L}$
there exists a line $L'\in{\mathcal L}'$ such that
$f(L)$ is contained in $L'$.
It is trivial that if $f$ and $f^{-1}$ both are semicollineations
then $f$ is a collineation.

\begin{theorem}[{\rm A. Kreuzer \cite{Kreuzer1}}]
If $\Pi$ satisfies the exchange axiom and
\begin{equation}
\dim\Pi \le\dim\Pi'
\end{equation}
then any semicollineation of $\Pi$ to $\Pi'$ is a collineation.
\end{theorem}

We will exploit the following lemma.

\begin{lemma}
Let $f:P\to P'$ be a collinearity preserving injection.
Then
$$f(\overline{X})\subset\overline{f(X)}$$
for any set $X\subset P$.
\end{lemma}

\begin{proof}
By Lemma 1.1, we need to show that
$$f([X]_{i})\subset [f(X)]_{i}$$
for each number $i\ge 0$.
We will prove this inclusion by induction.
For the case when $i=0$ it is trivial.
We suppose that $i\ge 1$ and take a point $p$ belonging to $[X]_{i}$.
Let $q$ and $q'$ be distinct points of $[X]_{i-1}$
such that $p$ lies on the line $qq'$.
By the inductive hypothesis,
$f(q)$ and $f(q')$ belong to $[f(X)]_{i-1}$.
Therefore,
$$f(p)\in f(q)f(q')\subset [f(X)]_{i}$$
and we get the required inclusion.
\end{proof}

\begin{proof}[Proof of Theorem 1.2]
We will prove that a semicollineation $f:P \to P'$ is non-collinearity preserving.
Let $p_{1},p_{2},p_{3}\in P$ be non-collinear points.
By Theorem 1.1, these points are contained in certain base
$B$ for $\Pi$.
Lemma 1.4 says that
$$P'=f(P)=f(\overline{B})\subset \overline{f(B)}.$$
Hence the linear space $\Pi'$ is spanned by the set $f(B)$
and (1.1) guarantees that $f(B)$ is a base for $\Pi'$;
in particular, $f(p_{1}),f(p_{2}),f(p_{3})$ are non-collinear.
\end{proof}

The following example shows that
the condition (1.1) in Theorem 1.2 can not be dropped
even if our linear spaces both are satisfying the exchange axiom.

\begin{exmp}[{\rm A. Kreuzer \cite{Kreuzer1}}]
{\rm
Suppose that $\Pi$ is a $3$-dimensi\-onal projective space
and $S$ is a plane of $\Pi$.
Define
$${\mathcal L}':=\{\;L\in{\mathcal L}\;|\;L\not\subset S\}\cup\{S\}.$$
The linear space $(P,{\mathcal L}')$ is spanned by the line $S$
and any point belonging to $P-S$; hence it is a plane.
It is not difficult to prove that this plane satisfies the exchange axiom.
The identical transformation of the set $P$ induces a semicollineation
of $\Pi$ to $(P,{\mathcal L}')$
which is not a collineation (since the dimensions of our spaces are different).
}\end{exmp}

It was noted above that any collineation maps bases to bases.
Now we show that collineations of linear spaces satisfying the exchange axiom
can be characterized as base preserving surjections.

\begin{theorem}[{\rm W.-l. Huang, A Krauzer \cite{HuangKreuzer}}]
Suppose that our spaces have the same dimension and satisfy the exchange axiom.
If $f:P\to P'$ is a surjection transferring each base for $\Pi$ to a base for $\Pi'$
then $f$ is a collineation of $\Pi$ to $\Pi'$.
\end{theorem}

\begin{proof}
Let us take two distinct points $p_{1},p_{2}\in P$
and a base $B$ for $\Pi$ containing them
(by Theorem 1.2, such base exists).
Then $f(B)$ is a base for $\Pi$.
Since our spaces have the same dimension, $f(p_{1})\ne f(p_{2})$.
Thus $f$ is injective.

The mapping $f$  is non-collinearity preserving;
indeed, three points of $P$ or $P'$
are non-collinear if and only if there exists a base for
the corresponding linear space containing them.
This implies that the inverse mapping $f^{-1}$ is collinearity preserving.
By Theorem 1.2, it is a collineation.
\end{proof}

\begin{rem}{\rm
Here Theorem 1.3 is given as a simple consequence of Theorem 1.2.
It was proved directly in \cite{HuangKreuzer}.
}\end{rem}

An injective mapping $f:P\to P'$ is called an {\it embedding}
if it is collinearity and non-collinearity preserving;
for this case $f$ transfers lines to subsets of lines and for any line
$L'\in {\mathcal L}'$ there is at most one line
$L\in {\mathcal L}$ such that $f(L)$ is contained in $L'$.
It is trivial that any bijective embedding is a collineation.

An embedding is said to be {\it strong} if it
maps independent subsets to independent subsets.
In the next subsection we give examples of non-strong embeddings.

\begin{rem}{\rm
Suppose that our linear spaces have the same dimension and
satisfy the exchange axiom.
Then strong embeddings of $\Pi$ to $\Pi'$ (if they exist) map bases to bases.
Inversely, let $f:P\to P'$ be a mapping which
sends each base for $\Pi$ to a base for $\Pi'$;
is this mapping a strong embedding? we do not know.
}\end{rem}

\subsection{Semilinear mappings and morphisms of projective spaces}
Let $V$ and $V'$ be left vector spaces over division rings $R$ and $R'$, respectively.
We require that the dimensions of these spaces are not less than $3$.
Then the sets $P(V)$ and $P(V')$ have the natural structure of
projective spaces (Example 1.4); these spaces are denoted by $\Pi(V)$ and $\Pi(V')$.

A mapping $l:V\to V'$ is said to be {\it semilinear} if
$$l(x+y)=l(x)+l(y)$$
for all $x,y\in V$ and there exists a homomorphism $\sigma: R\to R'$ such that
$$l(ax)=\sigma(a)l(x)$$
for all $x\in V$ and $a\in R$.
This homomorphism is uniquely defined if the mapping $l$ is non-zero;
for this case our mapping will be called also $\sigma$-{\it linear}.

Note that each non-zero homomorphism of a division ring is a
monomorphism and
for any monomorphism $\sigma:R\to R'$ there are non-zero
$\sigma$-linear mappings of $V$ to $V'$.
If $\sigma$ is an isomorphism then any $\sigma$-linear
mapping transfers subspaces to subspaces.
For semilinear mappings associated with non-surjective homomorphisms
this fails.

\begin{exmp}{\rm
The group ${\rm Aut}({\mathbb R})$ consists only of the identical transformation;
moreover, there are not non-surjective monomorphisms of ${\mathbb R}$ to itself
\cite{Benz}.
This means that any semilinear mapping between vector spaces over ${\mathbb R}$
is linear.
}\end{exmp}

\begin{exmp}{\rm
The group ${\rm Aut}({\mathbb C})$ is not trivial.
The transformation sending each complex number $z=a+bi$ to
the conjugate number $\overline{z}=a-bi$ is an automorphism;
but it is not unique non-identical automorphisms of the field ${\mathbb C}$.
Non-surjective monomorphisms of ${\mathbb C}$ to itself exist \cite{Benz}.
}\end{exmp}

\begin{exmp}{\rm
Let $F$ be a finite field.
Then $F$ consists of $p^{m}$ elements, where $p$ is a simple number.
For this case ${\rm Aut}(F)$ is a cyclic group of order $m$.
It is trivial that any monomorphism of $F$ to itself is surjective.
}\end{exmp}

We say that a $\sigma$-linear mapping is a {\it semilinear isomorphism} if
it is bijective and $\sigma$ is an isomorphism.
For any semilinear isomorphism the inverse mapping is a semilinear isomorphism.

It must be pointed out that there exist semilinear bijections
which are not semilinear isomorphisms.

\begin{exmp}{\rm
The semilinear bijection $l:{\mathbb R}^{2n} \to {\mathbb C}^{n}$
defined by the formula
$$l(x_{1},y_{1},\dots,x_{n},y_{n}):=(x_{1}+y_{1}i,\dots,x_{n}+y_{n}i)$$
is not a semilinear isomorphism.
}\end{exmp}

Let $l:V\to V'$ be a semilinear mapping.
The subspace consisting of all vectors $x\in V$ such that $l(x)=0$
is called the {\it kernel} of $l$ and denoted by ${\rm Ker}\,l$.
Our mapping is injective if and only if it's kernel is zero.
Now consider the mapping
$$P(l):P(V)-P({\rm Ker}\,l)\to P(V')$$
which transfers each $1$-dimensional subspace $Rx$, $x\in V-{\rm Ker}\,l$
to the $1$-dimensional subspace $R'l(x)$.
This mapping is globally defined only for the case when $l$ is injective.
However, the injectivity of $l$ does not guarantee that $P(l)$
is injective (Example 1.11).

\begin{rem}{\rm
Clearly, the mapping $P(l)$ is bijective if $l$ is a semilinear isomorphism.
Inversely, it is not difficult to prove that if $\dim V \le \dim V'$ and
$l:V\to V'$ is a semilinear mapping such that $P(l)$ is bijective
then $l$ is a semilinear isomorphism.
Is there for the case when $\dim V > \dim V'$
a semilinear mapping $l:V\to V'$ such that $P(l)$ is bijective?
}\end{rem}

For any semilinear isomorphism $l:V\to V'$ the mapping $P(l)$ is a collinea\-tion
of $\Pi(V)$ to $\Pi(V')$.
The classical version of the Fundamental Theorem of Projective Geometry
 says that
\begin{enumerate}
\item[]{\it if $\dim V\le \dim V'$ then
any semicollineation of $\Pi(V)$ to $\Pi(V')$
is induced by a semilinear isomorphism of $V$ to $V'$}
\end{enumerate}
(see \cite{Artin} or \cite{Baer}).
If the projective spaces $\Pi(V)$ and $\Pi(V')$ are isomorphic
then the vector spaces $V$ and $V'$ have the same dimension.
Thus
\begin{enumerate}
\item[]{\it all collineations of $\Pi(V)$ to $\Pi(V')$
are induced by semilinear isomorphisms of $V$ to $V'$}.
\end{enumerate}
This theorem was first proved by O. Veblen \cite{Veblen}
for the projective spaces over finite fields;
more historical information can be found in \cite{KarzelKroll}.

If $l:V\to V'$ is a semilinear injection then $P(l)$ is collinearity preserving
and globally defined.
Inversely, we have the following statement which is a consequence of
a more general result established by
C. A. Faure, A. Fr\"{o}lisher \cite{FaureFrolisher} and H. Havlicek \cite{Havlicek2}
(independently).

\begin{theorem}
Any colline\-arity preserving mapping of $P(V)$ to $P(V')$
is induced by a semilinear injection of $V$ to $V'$.
\end{theorem}

\begin{rem}{\rm
C. A. Faure, A. Fr\"{o}lisher and H. Havlicek have characterized partially defined
mappings induced by semilinear mappings.
A short proof of this statement can be found in Faure's paper \cite{Faure}
dedicated to Alfred Fr\"{o}lisher.
Collinearity preserving transformations of projective spaces were studied by many authors
(see, for example, \cite{Brauner1}, \cite{Sorensen}, \cite{CarterVogt});
more information related with geometrical characterizations of semilinear mappings
can be found in Havlicek's survey \cite{Havlicek4}.
}\end{rem}

Theorem 1.4 shows that
\begin{enumerate}
\item[]{\it any strong embedding of $\Pi(V)$ to $\Pi(V')$ is induced
by a semilinear injection of $V$ to $V'$ preserving the linear independence.}
\end{enumerate}
Non-strong embeddings exist.

\begin{exmp}{\rm (A. Brezuleanu, D. C. R\u{a}dulescu \cite{BrezuleanuRadulescu})
Let $F$ be a subfield of a field $F'$.
Assume also that there exist $b_{1},b_{2},b_{3}\in F'$ such that
$1,b_{1},b_{2},b_{3}$ are linearly independent vectors of the
vector space $F'$ over $F$.
The semilinear mapping $l:F^{4}\to F'^{3}$ defined by the formula
$$l(a_{1},a_{2},a_{3},a_{4})=(a_{1}+a_{4}b_{1},a_{2}+a_{4}b_{2},a_{3}+a_{4}b_{3})$$
is injective and sends each triple of linearly
independent vectors to linearly independent vectors.
This means that $P(l)$ is an embedding of $\Pi(F^{4})$ to $\Pi(F'^{3})$.
Clearly, this embedding is not strong.
}\end{exmp}

P. V. Ceccherini \cite{Ceccherini} has given an example of
a semicollineation of a $4$-dimensional projective space
to a non-Desargue\-sian projective plane.
On the other hand, any semicollineation of $\Pi(V)$ to $\Pi(V')$ is a collineation
if $\dim V \le \dim V'$.
Are there for the case when $\dim V>\dim V'$
semicollinea\-ti\-ons of $\Pi(V)$ to $\Pi(V')$?
Theorem 1.4 shows that this problem is equivalent to the
problem considered in Remark 1.4.

\newpage

\section{Grassmann spaces}

\setcounter{prop}{0}
\setcounter{equation}{0}
\setcounter{theorem}{0}
\setcounter{lemma}{0}
\setcounter{rem}{0}
\setcounter{exmp}{0}

In the course of this section we will suppose that
$\Pi=(P,{\mathcal L})$ and $\Pi'=(P',{\mathcal L}')$
are $n$-dimensional linear spaces satisfying the exchange axiom.

For each number $k\in \{0,1,\dots, n-1\}$ we denote by ${\mathcal G}_{k}(\Pi)$
the {\it Grassmann space} consisting of all $k$-dimensional subspaces of $\Pi$.
Then ${\mathcal G}_{0}(\Pi)=P$ and ${\mathcal G}_{1}(\Pi)={\mathcal L}$.
By Theorem 1.1, for any two elements $S$ and $U$ of
${\mathcal G}_{k}(\Pi)$ the inclusion $S\subset U$ implies that $S$ and $U$
are coincident.

\subsection{Adjacency relation}
Two elements of ${\mathcal G}_{k}(\Pi)$
are said to be {\it adjacent} if their intersection is $(k-1)$-dimensional.
By this definition, any two points ($0$-dimensional subspaces) are adjacent.
If $\Pi$ is projective then any two distinct $(n-1)$-dimensional subspaces are adjacent.

\begin{lemma}
If $U$ and $S$ are adjacent elements of ${\mathcal G}_{k}(\Pi)$ then
the subspace spanned by them is $(k+1)$-dimensional.
\end{lemma}

\begin{proof}
Let $B$ be a base for $S\cap U$.
Theorem 1.1 implies the existence of points $p\in S$ and $q\in U$ such that
$B\cup \{p\}$ and $B\cup \{q\}$ are bases for $S$ and $U$.
Then $B\cup \{p,q\}$ is a base for $\overline{S\cup U}$.
\end{proof}

Inversely, if our space is projective and
two $k$-dimensional subspaces span a $(k+1)$-dimensional subspace
then they are adjacent.
For the general case this fails:
there exists a plane satisfying the exchange axiom
and containing non-intersecting lines (Example 1.6).

\begin{rem}{\rm (A characterization of the adjacency relation in terms
of complements).
Two subspaces $S$ and $U$ of $\Pi$ are said to be {\it complementary}
if they span $\Pi$ and their intersection is empty;
for this case the dimension of $U$ is equal to $n-\dim S -1$.
If $\Pi$ is projective then for any two $k$-dimensional subspaces $S_{1}$ and $S_{2}$
the following conditions are equivalent:
\begin{enumerate}
\item[(1)] $S_{1}$ and $S_{2}$ are adjacent,
\item[(2)] there exists a subspace $S\in {\mathcal G}_{k}(\Pi)-\{S_{1},S_{2}\}$
such that each complement to $S$ is complementary to $S_{1}$ or $S_{2}$,
\end{enumerate}
(see A. Blunk, H. Havlicek \cite{BlunckHavlicek} and
H. Havlicek, M. Pankov \cite{HavlicekPankov}).
}\end{rem}

\begin{prop}
Any two $k$-dimensional subspaces $S$ and $U$ of $\Pi$ can be
connected by a sequence of adjacent $k$-dimensional subspaces;
in other words, there exists a sequence
\begin{equation}
S=S_{0},S_{1},\dots,S_{i}=U
\end{equation}
of elements of ${\mathcal G}_{k}(\Pi)$ where
$S_{j-1}$ and $S_{j}$ are adjacent for each number $j\in \{1,\dots,i\}$.
\end{prop}

\begin{proof}
Let us define
$$i:=k-\dim(S\cap U)$$
and prove by induction on $i$ the existence of
a sequence (2.1) of $k$-dimensional subspaces
such that $S_{j-1}$ and $S_{j}$ are adjacent for all $j\in \{1,\dots,i\}$.
If $i=1$ then $S$ and $U$ are adjacent and our statement is trivial.
For the general case any base $B$ for $S\cap U$ can be extended
to some bases $B_{S}$ and $B_{U}$ for $S$ and $U$ (Theorem 1.1).
We take two points
$$p\in B_{S}-B,\;q\in B_{U}-B$$
and consider the $k$-dimensional subspace
$$S_{1}:=\overline{(B_{S}-\{p\})\cup \{q\}}.$$
Then $S$ and $S_{1}$ are adjacent and
$$k-\dim (S_{1}\cap U)=i-1.$$
The inductive hypothesis implies the existence of a sequence $S_{2},\dots, S_{i}$
of $k$-dimensional subspaces satisfying the required conditions.
\end{proof}

\begin{rem}{\rm (The distance between subspaces)
Let $S$ and $U$ be distinct $k$-dimensional subspaces of $\Pi$.
The {\it distance} $d(S,U)$ between $S$ and $U$ can be defined as
the smallest number $i$ such that there exists a sequence (2.1)
satisfying the conditions of Proposition 2.1.
Then
$$d(S,U)\le k-\dim(S\cap U)$$
(see the proof of Proposition 2.1).
Now we show that
$$\dim\overline{S\cup U}-k\le d(S,U).$$
Indeed, suppose that (2.1) is a sequence of $k$-dimensional subspaces
where $S_{j-1}$ and $S_{j}$ are adjacent for all $j\in \{1,\dots,i\}$.
Then
$$\dim \overline{S_{0}\cup S_{1}}=k+1\;\mbox{ and }\;
\dim \overline{S_{0}\cup S_{j}}\le \dim \overline{S_{0}\cup S}_{j-1}+1.$$
Thus the dimension of $\overline{S\cup U}$ is not greater than $k+i$.
This implies the required inequality and we have established that
$$\dim\overline{S\cup U}-k\le d(S,U)\le k-\dim(S\cap U).$$
If our space is projective then
$$\dim \overline{S\cup U}=2k -\dim (S\cap U)$$
and we get the equality
$$d(S,U)=\dim\overline{S\cup U}-k=k-\dim(S\cap U)$$
which does not hold for the general case.
}\end{rem}

Let $U$ be a subspace of $\Pi$.
Denote by ${\mathcal G}_{k}(U)$ the set of
all $k$-dimensional subspaces incident to $U$;
recall that two subspaces are called {\it incident}
if one of these subspaces is contained in other.
If the subspace $U$ is not proper then ${\mathcal G}_{k}(U)$ coincides
with our Grassmann space.

Suppose that $U$ is $(k-1)$-dimensional.
Then any two distinct elements of the set ${\mathcal G}_{k}(U)$ are adjacent.
Moreover, if $k<n-1$ then for any $k$-dimensional subspace $S\notin {\mathcal G}_{k}(U)$
there exists $S'\in {\mathcal G}_{k}(U)$ such that
$S$ and $S'$ are not adjacent;
this means that ${\mathcal G}_{k}(U)$ is a maximal set of adjacent subspaces.
For the case when our space is projective and $k=n-1$
any two elements of ${\mathcal G}_{k}(\Pi)$
are adjacent and our set is not maximal.

If $\Pi$ is projective and $S$ is a $(k+1)$-dimensional subspace
then any two distinct elements of ${\mathcal G}_{k}(S)$ are adjacent;
moreover, ${\mathcal G}_{k}(S)$ is a maximal set of adjacent subspaces if $k>0$;
for $k=0$ this fails (since for this case any two distinct elements are adjacent).

\begin{prop}
Let ${\mathcal X}$ be a subset of ${\mathcal G}_{k}(\Pi)$
any two distinct element of which are adjacent.
Then there exists a $(k\pm 1)$-dimensional subspace $U$ such that
$${\mathcal X}\subset {\mathcal G}_{k}(U).$$
If $\Pi$ is projective then
any maximal set of adjacent subspaces coincides with
certain ${\mathcal G}_{k}(U)$ where $U$ is a $(k\pm 1)$-dimensional subspace of $\Pi$.
\end{prop}

The second assertion of Proposition 2.2 is old-known,
see \cite{Bertini}.
The following example shows that the same does not hold
if our space is not projective.

\begin{exmp}{\rm
Suppose that $\Pi$ has a plane $S$ containing two non-intersecting lines $L$ and $L'$
(see Example 1.6).
Let us take three distinct points $p_{1},p_{2}\in L$ and $q\in S-L$.
The lines
$$L,\;p_{1}q,\;p_{2}q$$
are mutually adjacent and Zorn Lemma implies the existence of
a maximal set of adjacent lines containing them.
This maximal set is a proper subset of ${\mathcal G}_{1}(S)$
(since the lines $L$ and $L'$ are not adjacent).
}\end{exmp}

\begin{proof}[Proof of Proposition 2.2]
Since the case $k=0,n-1$ is trivial, we will suppose that $0<k<n-1$.
Let $S_{1}$ and $S_{2}$ be distinct elements of the set ${\mathcal X}$.
Then
$$U:=\overline{S_{1}\cup S_{2}}\; \mbox{ and }\; U':=S_{1}\cap S_{2}$$
belong to ${\mathcal G}_{k+1}(\Pi)$ and ${\mathcal G}_{k-1}(\Pi)$, respectively.
If a $k$-dimensional subspace $S$ is adjacent to both $S_{1}$ and $S_{2}$
then it is incident to $U$ or $U'$.
Indeed, if $S$ is not contained in $U$ then the dimension of $S\cap U$ is less than $k$;
since $S\cap U$ contains the $(k-1)$-dimensional subspaces $S\cap S_{1}$ and $S\cap S_{2}$,
we have the equality
$$S\cap U=S\cap S_{1}=S\cap S_{2}$$
showing that $S\cap U$ coincides with $S_{1}\cap S_{2}=U'$.
Thus
$${\mathcal X}\subset {\mathcal G}_{k}(U)\cup {\mathcal G}_{k}(U').$$
If $S\in {\mathcal G}_{k}(U)$ and $S'\in {\mathcal G}_{k}(U')$ are adjacent
then one of these subspaces belongs to the intersection of
${\mathcal G}_{k}(U)$ and ${\mathcal G}_{k}(U')$
(suppose that $S'$ does not belong to ${\mathcal G}_{k}(U)$,
then $U\cap S'=U'$ contains $S\cap S'$, these subspaces both are
$(k-1)$-dimensional and we have $S\cap S'=U'$).
Therefore, ${\mathcal X}$ is a subset of ${\mathcal G}_{k}(U)$ or ${\mathcal G}_{k}(U')$.

Now suppose that our space is projective.
Then ${\mathcal G}_{k}(U)$ and ${\mathcal G}_{k}(U')$ are maximal sets of adjacent subspaces.
If the same holds for ${\mathcal X}$
then ${\mathcal X}$ coincides with one of these sets.
\end{proof}

\begin{rem}{\rm
A. Beutelspacher, J. Eisfeld, J. M\"{u}ller
\cite{BeutelspacherEisfeldMuller}
have studied subsets of Grassmann spaces such that
the distance between any two distinct elements is equal to $2$.
}\end{rem}

\begin{rem}{\rm
If $0<k<n-1$ and $\Pi$ is projective
then ${\mathcal G}_{k}(\Pi)$ has the natural structure of a {\it partial linear space}
(there is a set of proper subsets called lines such that each point belongs
to certain line, each line contains at least two points,
but pairs of non-collinear points exist and for any two distinct points
there is at most one line containing them);
lines of this space are defined as
$${\mathcal G}_{k}(S)\cap{\mathcal G}_{k}(U)$$
where $U$ and $S$ are incident elements of ${\mathcal G}_{k-1}(\Pi)$
and ${\mathcal G}_{k+1}(\Pi)$.
Then two points are collinear if and only if they are
adjacent subspaces.
K. Pra\.{z}movski and M. \.{Z}ynel \cite{Prazmovski} and \cite{PrazmovskiZynel}
study so-called {\it spine spaces}:
for a fixed subspace $S$ of $\Pi$ and some fixed number
$m\in\{-1,\dots,\dim S\}$
consider the set ${\mathcal F}_{km}(S)$ consisting of
all subspaces $U\in {\mathcal G}_{k}(\Pi)$
such that
$$\dim U\cap S=m,$$
a subset of ${\mathcal F}_{km}(S)$ is line
if it contains at least two elements and is the intersection
of some line of ${\mathcal G}_{k}(\Pi)$ with  ${\mathcal F}_{km}(S)$.
}\end{rem}

\subsection{Base subsets}
If $B$ is a base for $\Pi$ then
the set consisting of all $k$-dimensional subspaces spanned
by points of $B$ is said to be the {\it base} subset of ${\mathcal G}_{k}(\Pi)$
(associated with the base $B$ or defined by $B$).
Each base subset of ${\mathcal G}_{k}(\Pi)$ consists of $\binom{n+1}{k+1}$ elements.
It is trivial that any base subset of ${\mathcal G}_{0}(\Pi)=P$ is a base for $\Pi$.

\begin{prop}
Let $S$ and $U$ be distinct subspaces of $\Pi$.
The following conditions are equivalent:
\begin{enumerate}
\item[{\rm (1)}] there exists a base for $\Pi$ such that $S$ and $U$ are spanned
by points of this base,
\item[{\rm (2)}] $\dim \overline{S\cup U}=\dim S +\dim U -\dim (S\cap U)$.
\end{enumerate}
\end{prop}

\begin{proof}
The implication $(1)\Rightarrow (2)$ follows from Lemma 1.2.

$(2)\Rightarrow (1)$.
Let us take a base for $S\cap U$ and extend it to
bases $B_{S}$ and $B_{U}$ of $S$ and $U$, respectively.
The subspace $\overline{S\cup U}$ is spanned by the set $B_{S}\cup B_{U}$.
The condition (2) shows that this set is a base for $\overline{S\cup U}$.
By Theorem 1.1, this base can be extended to a base for $\Pi$.
\end{proof}

The following assertions are simple consequences of Proposition 2.3:
\begin{enumerate}
\item[{\rm (1)}]
If $\Pi$ is projective then
for any two subspaces there exists a base for $\Pi$ such that
these subspaces are spanned by points of this base.
\item[{\rm (2)}]
Two $k$-dimensional subspaces $S$ and $U$ of $\Pi$
are contained in certain base subset of ${\mathcal G}_{k}(\Pi)$
if and only if
$$\dim \overline{S\cup U}=2k -\dim (S\cap U);$$
in particular, for any two adjacent elements of ${\mathcal G}_{k}(\Pi)$
there exists a base subset containing them.
\end{enumerate}

Let $V$ be an $(n+1)$-dimensional vector space over a field $F$.
Then $\Pi(V)$ is an $n$-dimensional pappian projective space.
Consider the exterior power vector space $\Lambda^{k+1}(V)$ and
the associated $(\binom{n+1}{k+1}-1)$-dimensional projective space
$\Pi(\Lambda^{k+1}(V))$.
Any base
$$Fx_{1},\dots, Fx_{n+1}$$
for $\Pi(V)$ defines the base for $\Pi(\Lambda^{k+1}(V))$
consisting of all
$$F(x_{i_{1}}\wedge\dots\wedge x_{i_{k+1}}).$$
Clearly, there exist bases for $\Pi(\Lambda^{k+1}(V))$ which
are not associated with any base of $\Pi(V)$.
The Grassmann injection \cite{HodgePedoe}
$${\mathcal G}_{k}(\Pi(V))\to P(\Lambda^{k+1}(V))$$
$$Fy_{1}+\dots+Fy_{k+1}\to F(y_{1}\wedge\dots\wedge y_{k+1})$$
sends the base subset associated with some base for $\Pi(V)$
to the base for $\Pi(\Lambda^{k+1}(V))$ defined by this base.

\subsection{Principles of duality for projective spaces}
Now suppose that $\Pi$ is projective and
define $P^{*}:={\mathcal G}_{n-1}(\Pi)$.
A subset of $P^{*}$ is said to be a {\it line} if it consists of all
$p^{*}\in P^{*}$ containing certain $(n-2)$-dimensional subspace of $\Pi$.
Then the following statements hold true:
\begin{enumerate}
\item[(1)]
The set $P^{*}$ together with the family of lines defined above
is an $n$-dimen\-si\-onal projective space.
This projective space is called {\it dual} to $\Pi$ and denoted by $\Pi^{*}$.
\item[(2)]
A subset of $P^{*}$ is a $k$-dimensional subspace of the dual space if
and only if it consists of all
$p^{*}\in P^{*}$ containing certain $(n-k-1)$-dimensional subspace of $\Pi$.
Thus there is the natural bijection between
${\mathcal G}_{k}(\Pi^{*})$ and ${\mathcal G}_{n-k-1}(\Pi)$.

\item[(3)]
A subset of $P^{*}$ is a base for the dual space
if and only if it consists of all $(n-1)$-dimensional subspaces
spanned by points of certain base for $\Pi$,
in other words, it is a base subset of ${\mathcal G}_{n-1}(\Pi)$.
If $B$ is a base for $\Pi$ and $B^{*}$ is the base for $\Pi^{*}$ defined by $B$
then the base subset of ${\mathcal G}_{n-k-1}(\Pi)$ associated with $B$ is
the base subset of ${\mathcal G}_{k}(\Pi^{*})$ associated with $B^{*}$.
\end{enumerate}
The natural bijection of ${\mathcal G}_{k}(\Pi^{*})$ to ${\mathcal G}_{n-k-1}(\Pi)$
(see (2)) is distance preserving;
by (3), it preserves the class of base subsets.
In what follows we will suppose that these Grassmann spaces are coincident;
in particular, the second dual space $\Pi^{**}$ coincides with $\Pi$.

\begin{exmp}{\rm
Let $V$ be an $n$-dimensional vector space over a division ring and $n\ge 3$.
The projective space $\Pi(V)^{*}$ is canonically isomorphic to $\Pi(V^{*})$
where $V^{*}$ is the vector space dual to $V$.
Indeed, points of $\Pi(V)^{*}$
can be considered as $(n-1)$-dimensional subspaces of $V$,
the bijection sending each $(n-1)$-dimensional subspace of $V$
to its annihilator defines a collineation of $\Pi(V)^{*}$  to $\Pi(V^{*})$.
}\end{exmp}

\subsection{Mappings induced by strong embeddings}

Let $f$ be a strong embedding of $\Pi$ to $\Pi'$.
Then for any subspace $S$ of $\Pi$ we have
$$\dim S= \dim \overline{f(S)}.$$
Thus $f$ induces the mapping
$$G_{k}(f):{\mathcal G}_{k}(\Pi)\to {\mathcal G}_{k}(\Pi')$$
$$S\to\overline{f(S)}.$$
It is easy to see that this mapping is injective.

\begin{prop}
If $\Pi$ is projective
and for certain number $k$ the mapping $G_{k}(f)$ is bijective then $f$ is a collineation.
\end{prop}

\begin{proof}
Let $U'$ be an arbitrary element of ${\mathcal G}_{k-1}(\Pi')$.
We take two $k$-dimensi\-on\-al subspaces $S'_{1}$ and $S'_{2}$ of $\Pi'$ such that
$$U'=S'_{1}\cap S'_{2}$$
(the existence of these subspaces follows from Theorem 1.1).
Then $S'_{1}$ and $S'_{2}$ are adjacent and
the subspace spanned by them is $(k+1)$-dimensional.
The mapping $G_{k}(f)$ is bijective and the equalities
$$\overline{f(S_{1})}=S'_{1}\;\mbox{ and }\;\overline{f(S_{2})}=S'_{2}$$
hold for some $k$-dimensional subspaces $S_{1}$ and $S_{2}$ of $\Pi$.
Lemma 1.4 gives the inclusion
$$f(\overline{S_{1}\cup S_{2}})\subset \overline{f(S_{1})\cup f(S_{2})}\subset
\overline{S'_{1}\cup S'_{2}}$$
showing that
the dimension of $\overline{S_{1}\cup S}_{2}$ is not greater than $k+1$.
Since $S_{1}$ and $S_{2}$ are distinct $k$-dimensional subspaces,
this dimension is equal to $k+1$.
Hence $S_{1}$ and $S_{2}$ are adjacent (the space $\Pi$ is projective)
and
$$U:=S_{1}\cap S_{2}$$
is an element of ${\mathcal G}_{k-1}(\Pi)$.
Then
$$\overline{f(U)}=\overline{f(S_{1})\cap f(S_{2})}\subset
\overline{f(S_{1})}\cap \overline{f(S_{2})}=S'_{1}\cap S'_{2}=U'$$
(the intersection of $\overline{f(S_{1})}$ and $\overline{f(S_{2})}$
is a subspace containing $f(S_{1})\cap f(S_{2})$, thus
$\overline{f(S_{1})\cap f(S_{2})}$ is contained in
$\overline{f(S_{1})}\cap \overline{f(S_{2})}$\,).
The subspaces $\overline{f(U)}$ and $U'$ both are $(k-1)$-dimensional
and
$$\overline{f(U)}=U'.$$
We have established that $G_{k-1}(f)$ is bijective.
Step by step we can prove that $G_{0}(f)=f$ is bijective.
\end{proof}

The following example shows that the statement given above fails
if $\Pi$ is not projective.

\begin{exmp}{\rm
Suppose that $\Pi'$ is a projective space and $\Pi=\Pi'_{X}$
where $X=P'-\{p\}$ (see Example 1.6).
The mapping of ${\mathcal G}_{1}(\Pi)={\mathcal L}'_{X}$ to
${\mathcal G}_{1}(\Pi')={\mathcal L}'$
induced by the natural embedding of $\Pi'_{X}$ to $\Pi'$ is bijective.
}\end{exmp}

If $\Pi'$ is projective and $g$ is a strong embedding of $\Pi$ to $\Pi'^{*}$
then, by the principle of duality,
$$G_{k}(g):{\mathcal G}_{k}(\Pi)\to {\mathcal G}_{k}(\Pi'^{*})$$
can be considered as an injection to
${\mathcal G}_{n-k-1}(\Pi')$; for the case when $n=2k+1$ we get an injective mapping
of ${\mathcal G}_{k}(\Pi)$ to ${\mathcal G}_{k}(\Pi')$.

Now suppose that our spaces both are projective.
Then for any strong embedding $f$ of $\Pi$ to $\Pi'$
the mapping $G_{n-1}(f)$ is a strong embedding of $\Pi^{*}$ to $\Pi'^{*}$,
this mapping will be called the {\it contragradient} to $f$ and denoted by ${\check f}$.
The principle of duality says that
$$G_{k}({\check f})=G_{n-k-1}(f).$$
By Proposition 2.4, ${\check f}$ is a collineation if
and only if $f$ is a collineation.

\begin{rem}{\rm
Let $V$ and $V'$ be left vector spaces over isomorphic division rings
and $l:V\to V'$ be a $\sigma$-linear isomorphism.
Recall that the {\it contragradient} of $l$ is the semilinear isomorphism
${\check l}:V^{*}\to V'^{*}$ sending
any linear functional $\alpha$ to the functional defined by
the formula
$$x'\to \sigma(\alpha(l^{-1}(x')))\;\;\;\;\;\forall\;x'\in V'.$$
Now suppose that the dimensions of $V$ and $V'$ are not less than $3$
and consider the collineation $f$ of $\Pi(V)$ to $\Pi(V')$ induced by $l$.
Since $\Pi(V)^{*}$ and $\Pi(V')^{*}$ are canonically isomorphic to
$\Pi(V^{*})$ and $\Pi(V'^{*})$, respectively (see Example 2.2),
${\check f}$ defines certain collineation of $\Pi(V^{*})$ to $\Pi(V'^{*})$.
This collineation is induced by ${\check l}$.
}\end{rem}

\subsection{Adjacency preserving transformations}
In this section we will always suppose that the spaces $\Pi$ and $\Pi'$ are projective.
Now we give a few examples of distance preserving bijections
of ${\mathcal G}_{k}(\Pi)$ to ${\mathcal G}_{k}(\Pi')$.
By the definition of the distance between subspaces (Remark 2.2),
a bijection of ${\mathcal G}_{k}(\Pi)$ to ${\mathcal G}_{k}(\Pi')$
is distance preserving if and only if
it is adjacency preserving in both directions.
\begin{enumerate}
\item[(1)]
If $k=0,n-1$ then any bijection of ${\mathcal G}_{k}(\Pi)$ to ${\mathcal G}_{k}(\Pi')$
is distance preserving.
\item[(2)]
For any collineation  $f$ of $\Pi$ to $\Pi'$
the mapping $G_{k}(f)$ preserves the distant between subspaces.
\item[(3)]
Recall that he canonical bijection of
${\mathcal G}_{k}(\Pi)$ to ${\mathcal G}_{n-k-1}(\Pi^{*})$
(see Subsection 2.3) is distance preserving.
Thus if $n=2k+1$ then any bijection of
${\mathcal G}_{k}(\Pi)$ to ${\mathcal G}_{k}(\Pi')$
induced by a collineation of $\Pi$ to $\Pi'^{*}$
is distance preserving.
\end{enumerate}
Chow's Theorem \cite{Chow} says that there are not other
distance preserving bijections.

\begin{theorem}[{\rm W.L. Chow \cite{Chow}}]
Let $0<k<n-1$ and $f$ be a distance preserving bijection of ${\mathcal G}_{k}(\Pi)$
to ${\mathcal G}_{k}(\Pi')$.
Then the following statements hold true:
\begin{enumerate}
\item[{\rm (1)}]
If $n\ne 2k+1$ then $f$ is induced by a collineation of $\Pi$ to $\Pi'$.
\item[{\rm (2)}]
If $n=2k+1$ then $f$ is induced by a collineation of $\Pi$ to $\Pi'$
or $\Pi'^{*}$.
\end{enumerate}
\end{theorem}

\begin{proof}[Sketch of proof {\rm \cite{Chow}, \cite{Dieudonne}}]
The mapping $f$ preserves the class of maximal sets of adjacent subspaces.
By Proposition 2.2, there are two different types of maximal sets
(associated with $(k-1)$-dimensional and $(k+1)$-dimensional subspaces,
respectively).
The intersection of two maximal sets has only one element
if and only if these sets have the same type and
the corresponding subspaces are adjacent.
Since any two elements of Grassmann space can be connected by
a sequence of adjacent elements,
only one of the following possibilities is realized:
\begin{enumerate}
\item[(A)]
$f$ preserves type of any maximal set of adjacent subspaces,
\item[(B)]
$f$ changes types of all maximal sets.
\end{enumerate}
{\it Case} (A).
The mapping $f$ induces some distance preserving
bijection of ${\mathcal G}_{k-1}(\Pi)$ to ${\mathcal G}_{k-1}(\Pi')$
which is a collineation if $k=1$.
Thus we can prove by induction that for this case
$f$ is induced by a collineation of $\Pi$ to $\Pi'$.

{\it Case} (B).
We have a distance preserving bijection of
${\mathcal G}_{k-1}(\Pi)$ to ${\mathcal G}_{k+1}(\Pi')$;
it is possible only for $n=2k+1$.
Then $f$ can be considered as a distance preserving bijection of
${\mathcal G}_{k}(\Pi)$ to ${\mathcal G}_{k}(\Pi'^{*})$ which satisfies
the condition (A);
this means that $f$ is induced by a collineation of $\Pi$ to $\Pi'^{*}$.
\end{proof}

W.-l. Huang \cite{Huang1} has shown that
{\it each adjacency preserving bijection of ${\mathcal G}_{k}(\Pi)$ to
${\mathcal G}_{k}(\Pi')$ preserves the distance between subspaces};
therefore, {\it if $0<k<n-1$ then the statements {\rm (1)} and {\rm (2)} of
Theorem 3.1 are fulfilled for any adjacency preserving bijection
$f:{\mathcal G}_{k}(\Pi)\to {\mathcal G}_{k}(\Pi')$}.
For $k=1$ this statement
was proved by H. Brauner \cite{Brauner2} and H. Havlicek \cite{Havlicek3},
but their methods can not be used for the general case.
If $\Pi$ is not projective then this result does not hold;
for this case there is an adjacency preserving bijection of
${\mathcal G}_{1}(\Pi)$ to ${\mathcal G}_{1}(\Pi')$ such that
the inverse mapping is not adjacency preserving (Example 2.3).

A. Kreuzer \cite{Kreuzer2} has constructed an example
of a bijection of the set of lines of an $n$-dimensional projective space
($n\ge 3$) to the set of lines of a projective plane.
Since any two lines of a projective plane are adjacent,
this bijection is adjacency preserving; however, for the inverse bijection
this fails.
This plane is Desarguesian and can be considered as a subspace of
some $n$-dimensional projective space.
Therefore, {\it there exists an adjacency preserving injection
of the set of lines of an $n$-dimensional projective space to
the set of lines of another $n$-dimensional projective space
which is not induced by a strong embedding.}

An analogue of Chow's Theorem for linear spaces
was given by H. Havlicek \cite{Havlicek5}.
He considered adjacency preserving bijective transformations of
the set of lines of a linear space and did not require that this space
is satisfying the exchange axiom.
If the dimension of the linear space is not less than $4$
then these mappings are induced by collineations;
for $3$-dimensional linear spaces Havlicek's result is not so simple-formulated
and we do not give it here.

\begin{rem}{\rm
If $0<k<n-1$ then the Grassmann spaces ${\mathcal G}_{k}(\Pi)$
has the natural structure of a partial linear space
such that two points of this spaces are collinear if and only if they are
adjacent subspaces (Remark  2.4).
Collinearity preserving mappings of Grassmann spaces to projective spaces
were studied by H. Havlicek \cite{Havlicek1}, A. L. Wells \cite{Wells} and
C. Zanella \cite{Zanella}.
}\end{rem}

\subsection{Transformations preserving base subsets}
The mappings of Grassmann spaces induced by strong embeddings
transfer base subsets to base subsets:
\begin{enumerate}
\item[(1)]
Let $f$ be a strong embedding of $\Pi$ to $\Pi'$.
If $B$ is a base for $\Pi$ then $f(B)$ is a base for $\Pi'$
and the mapping $G_{k}(f)$ sends the base subset associated with $B$ to
the base subset associated with $f(B)$.
\item[(2)]
Suppose that $\Pi'$ is projective and $n=2k+1$.
Then for any strong embedding $g$ of $\Pi$ to $\Pi'^{*}$ the mapping $G_{k}(g)$ is
an injection of ${\mathcal G}_{k}(\Pi)$ to ${\mathcal G}_{k}(\Pi')$
sending base subsets to base subsets
(indeed, if $B$ is a base for $\Pi$ then $g(B)$ is a base for $\Pi'^{*}$,
the base $g(B)$ is defined by certain base $B'$ for $\Pi'$
and $G_{k}(g)$ maps the base subset associated with $B$ to
the base subset associated with $B'$).
\end{enumerate}
We show that there are not other mappings satisfying this condition
if our spaces both are projective.

\begin{theorem}
Let $\Pi$ and $\Pi'$ be projective and $0<k<n-1$.
Let also $f$ be a mapping of ${\mathcal G}_{k}(\Pi)$ to ${\mathcal G}_{k}(\Pi')$
which sends base subsets to base subsets.
Then the following statements hold true:
\begin{enumerate}
\item[{\rm(1)}]
If $n\ne 2k+1$ then $f$ is induced by a strong embedding of $\Pi$ to $\Pi'$.
\item[{\rm(2)}]
If $n=2k+1$ then
$f$ is induced by a strong embedding of $\Pi$ to $\Pi'$ or $\Pi'^{*}$.
\end{enumerate}
\end{theorem}

\begin{rem}{\rm
Unfortunately, we can not prove the same statement for the case when
$k=0,n-1$ (see Remark 1.3).
}\end{rem}

If $\Pi$ and $\Pi'$ both are projective and $n<2k+1$ then
$$n>2(n-k-1)+1$$
and any mapping $f$ of ${\mathcal G}_{k}(\Pi)$ to ${\mathcal G}_{k}(\Pi')$
transferring base subsets to base subsets
can be considered as a mapping
of ${\mathcal G}_{n-k-1}(\Pi^{*})$ to ${\mathcal G}_{n-k-1}(\Pi'^{*})$
satisfying the same condition.
If the last mapping is induced by a strong embedding of $\Pi^{*}$ to $\Pi'^{*}$
then $f$ is induced by the contragradient of this embedding.
Therefore, we need to prove Theorem 2.2 only for the case when $n\ge 2k+1$.

\begin{cor}
Let $\Pi$ and $\Pi'$ be as in Theorem 2.2
and $f$ be a surjection of ${\mathcal G}_{k}(\Pi)$ to ${\mathcal G}_{k}(\Pi')$
sending base subsets to base subsets.
Then the following statements are fulfilled:
\begin{enumerate}
\item[{\rm(1)}]
If $n\ne 2k+1$ then $f$ is induced by a collineation of $\Pi$ to $\Pi'$.
\item[{\rm(2)}]
If $n=2k+1$ then $f$ is induced by a collineation of $\Pi$ to $\Pi'$ or $\Pi'^{*}$.
\end{enumerate}
\end{cor}

For the cases $k=0,n-1$ this is a partial case of Theorem 1.3.
If $1<k<n-1$ then this statement follows from Theorem 2.2 and Proposition 2.4
(it was proved directly in \cite{Pankov3}, some weak versions of this result
can be found in \cite{Pankov1} and \cite{Pankov2}).

\begin{theorem}
Let $\Pi$ and $\Pi'$ be satisfying the axiom
{\rm (P2)}\footnote{this axiom says that each line has at least three distinct points}
and $n\ge 4$.
If $f:{\mathcal L}\to{\mathcal L}'$ is an injection
transferring base subsets to base subsets then
$f$ is induced by a strong embedding of $\Pi$ to $\Pi'$.
\end{theorem}

Example 2.3 shows that bijections sending base subsets to base subsets
are not induced by collineation.
The bijection constructed in Example 2.3
maps base subsets to base subsets, but the inverse mapping
does not satisfy this condition.

\begin{theorem}
Let $\Pi$ and $\Pi'$ be as in Theorem 2.3 and $n>2k+1$.
Let also $f$ be a bijection of ${\mathcal G}_{k}(\Pi)$ to ${\mathcal G}_{k}(\Pi')$
such that $f$ and $f^{-1}$ send base subsets to base subsets.
Then $f$ is induced by a collineation of $\Pi$ to $\Pi'$.
\end{theorem}

We always require that our linear spaces have the same dimensional.
If $\dim \Pi < \dim \Pi'$ then there
exist injective mappings of ${\mathcal G}_{k}(\Pi)$ to ${\mathcal G}_{k}(\Pi')$
which transfer base subsets to subsets of base subsets and are not induced
by strong embeddings of $\Pi$ to $\Pi'$.

\begin{exmp}{\rm
Suppose that $n\ge 4$.
Let $S$ be an $(n-2)$-dimensional subspace of $\Pi$ and $L$ be a complement to $S$.
Clearly, $L$ is a line.
We take an arbitrary line $L'\subset S$
and consider the mapping of ${\mathcal G}_{1}(\Pi_{S})$ to ${\mathcal G}_{1}(\Pi)$
which is identical on the set ${\mathcal G}_{1}(\Pi_{S})-\{L'\}$ and
transfers $L'$ to $L$.
It is an injection sending any base subset of ${\mathcal G}_{1}(\Pi_{S})$
to a subset of a base subset of ${\mathcal G}_{1}(\Pi)$
(indeed, for any base subset ${\mathcal B}$ of ${\mathcal G}_{1}(\Pi_{S})$
there is a base subset of ${\mathcal G}_{1}(\Pi)$ containing  ${\mathcal B}\cup \{L\}$).
However, our mapping is not induced by an embedding of $\Pi_{S}$ to $\Pi$.
}\end{exmp}

The results given above (Theorems 2.2 - 2.4) will be proved in two steps.
First of all we establish the following.

\begin{theorem}
If $\Pi$ and $\Pi'$ satisfy the axiom {\rm (P2)} and $0<n<n-1$
then any injection of ${\mathcal G}_{k}(\Pi)$ to ${\mathcal G}_{k}(\Pi')$
sending base subsets to base subsets is adjacency preserving.
\end{theorem}

The second step is a modification of the proof of Theorem 2.1.

\subsection{Proof of Theorem 2.5}

In this subsection we establish a few properties of base subsets
and use them to prove Theorem 2.5.

Let $\Pi$ be satisfying (P2) and $0<k<n-1$.
Let also $B=\{p_{1},\dots,p_{n+1}\}$ be a base for $\Pi$.
For any number $m\in \{0,1,\dots, n-1\}$ we denote by ${\mathcal B}_{m}$
the base subset of ${\mathcal G}_{m}(\Pi)$ associated with $B$
(then ${\mathcal B}_{0}$ coincides with $B$).

Let us put ${\mathcal B}_{k}(+i)$ and ${\mathcal B}_{k}(-i)$
for the sets of all elements of ${\mathcal B}_{k}$ which contain $p_{i}$
or do not contain $p_{i}$, respectively.
More general, if $S\in {\mathcal B}_{m}$ then we denote by ${\mathcal B}_{k}(S)$ the set of all
elements of ${\mathcal B}_{k}$ incident to $S$;
for the case when $m=n-1$ there exists unique number $j$ such that $p_{j}\not\in S$
and ${\mathcal B}_{k}(S)$ coincides with ${\mathcal B}_{k}(-j)$.
A direct verification shows that
$$|{\mathcal B}_{k}(S)|=
\begin{cases}
\binom{m+1}{k+1}&\mbox{ if }\;m\ge k \\
\binom{n-m}{k-m}&\mbox{ if }\;m<k;
\end{cases}
$$
in particular,
$${\mathcal B}_{k}(+i)=\binom{n}{k}\;\mbox{ and }\;
{\mathcal B}_{k}(-i)=\binom{n}{k+1}.$$

Let ${\mathcal R}$ be a subset of ${\mathcal B}_{k}$.
We say that ${\mathcal R}$ is {\it exact} if ${\mathcal B}_{k}$ is
unique base subset of ${\mathcal G}_{k}(\Pi)$ containing ${\mathcal R}$;
otherwise, the set ${\mathcal R}$ is said to be {\it inexact}.
Denote by $S_{i}({\mathcal R})$
the intersection of all elements of ${\mathcal R}$ containing $p_{i}$.
It is easy to see that if
\begin{equation}
S_{i}({\mathcal R})=p_{i}
\end{equation}
for all numbers $i$ then ${\mathcal R}$ is exact.
The inverse statement follows from the axiom (P2).

\begin{lemma}
${\mathcal R}$ is an exact subset of ${\mathcal B}_{k}$
if and only if we have $(2.2)$ for each number $i$.
\end{lemma}

\begin{proof}
Suppose that (2.2) does not hold for certain $i$.
Then one of the following possibilities is realized:
\begin{enumerate}
\item[(A)] $S_{i}({\mathcal R})$ contains some point $p_{j}$ such that $j\ne i$,
\item[(B)] $S_{i}({\mathcal R})=\emptyset$.
\end{enumerate}
For the case (A) we choose a point
$$p\in p_{i}p_{j}-\{p_{i},p_{j}\}$$
(by the axiom (P2), this point exists).
For the second case we can take an arbitrary point $p\ne p_{i}$ belonging to
$P-\overline{B-\{p_{i}\}}$.
For each of these cases the base subset of ${\mathcal G}_{k}(\Pi)$
associated with the base
$$(B-\{p_{i}\})\cup \{p\}$$
contains ${\mathcal R}$;
this means that ${\mathcal R}$ is inexact.
\end{proof}

\begin{exmp}{\rm
Let us take two distinct numbers $i,j$ and
suppose that that
\begin{equation}
{\mathcal R}={\mathcal B}_{k}(-i)\cup {\mathcal B}_{k}(p_{i}p_{j}).
\end{equation}
Then $S_{l}({\mathcal R})=p_{l}$ if $l\ne i$ and
$S_{i}({\mathcal R})=p_{i}p_{j}$.
Thus our set is inexact.
Note that for any $U\in {\mathcal B}_{k}-{\mathcal R}$
the subset ${\mathcal R}\cup \{U\}$ is exact
(indeed, $U\cap S_{i}({\mathcal R})=p_{i}$).
This implies that ${\mathcal R}$ is a maximal inexact subset.
Since ${\mathcal B}_{k}(-i)$ and ${\mathcal B}_{k}(p_{i}p_{j})$
are non-intersecting, $|{\mathcal R}|$
is equal to $\binom{n}{k+1}+\binom{n-1}{k-1}$.
}\end{exmp}

\begin{lemma}
If ${\mathcal R}$ is a maximal inexact subset of ${\mathcal B}_{k}$
then there exist two distinct numbers $i$ and $j$
such that the equality $(2.3)$ holds true.
\end{lemma}

\begin{proof}
Since ${\mathcal R}$ is inexact,
for certain number $i$ the subspace $S_{i}({\mathcal R})$
is empty or contains at least two points of the base $B$.
If $S_{i}({\mathcal R})$ is not empty then
we choose $j\ne i$ such that $p_{j}\in S_{i}({\mathcal R})$;
for the case when  $S_{i}({\mathcal R})=\emptyset$ we can take any
number $j\ne i$.
Then for any $U\in {\mathcal R}$ one of the following possibilities is realized:
\begin{enumerate}
\item[(A)] $p_{i}\in U$ then $p_{i}p_{j}\subset S_{i}({\mathcal R})\subset U$ and
$U\in {\mathcal B}_{k}(p_{i}p_{j})$,
\item[(B)] $p_{i}\notin U$ then $U\in {\mathcal B}_{k}(-i)$.
\end{enumerate}
Hence
$${\mathcal R}\subset {\mathcal B}_{k}(-i)\cup {\mathcal B}_{k}(p_{i}p_{j}).$$
Since our inexact set is maximal,
we get the required equality.
\end{proof}

The arguments used to prove Lemma 2.3 also shows that
\begin{enumerate}
\item[]{\it an inexact subset of ${\mathcal B}_{k}$ is maximal if and only if it consists of
$$\binom{n}{k+1}+\binom{n-1}{k-1}$$ elements}.
\end{enumerate}

>From this moment the complements to maximal inexact subsets will be called {\it complement subsets}.
In other words,
${\mathcal R}\subset {\mathcal B}_{k}$ is said to be a {\it complement subset}
of ${\mathcal B}_{k}$ if ${\mathcal B}_{k}-{\mathcal R}$ is a maximal inexact set.
For this case Lemma 2.3 implies the existence of two distinct numbers  $i$ and $j$ such that
$${\mathcal B}_{k}-{\mathcal R}={\mathcal B}_{k}(-i)\cup {\mathcal B}_{k}(p_{i}p_{j}).$$
Then
$${\mathcal R}={\mathcal B}_{k}(+i)\cap{\mathcal B}_{k}(-j);$$
we denote this complement set by ${\mathcal B}_{k}(+i,-j)$.

We define
$$m:=\max\{k, n-k-1\}$$
and say that $m+1$ complement sets
$${\mathcal R}^{1}={\mathcal B}_{k}(+i_{1},-j_{1}),\dots,
{\mathcal R}^{m+1}={\mathcal B}_{k}(+i_{m+1},-j_{m+1})$$
form a {\it regular} collection if their intersection consists of one element.
Clearly, this intersection is not empty only if
\begin{equation}
\{i_{1},\dots,i_{m+1}\}\cap\{j_{1},\dots,j_{m+1}\}=\emptyset .
\end{equation}

\begin{lemma}
The collection ${\mathcal R}^{1},\dots, {\mathcal R}^{m+1}$
is regular if and only if $(2.4)$ holds and one of
the following possibilities is realized:
\begin{enumerate}
\item[]$n>2k+1$ (then $m=k$) and $i_{1},\dots,i_{k+1}$ are different,
\item[]$n=2k+1$ (then $m=k=n-k-1$) and $i_{1},\dots,i_{k+1}$ or $j_{1},\dots,j_{k+1}$
are different,
\item[]$n<2k+1$ (then $m=n-k-1$) and $j_{1},\dots,j_{n-k}$ are different.
\end{enumerate}
\end{lemma}

\begin{proof}
Direct verification.
\end{proof}

A collection of $m$ complement subsets of ${\mathcal B}_{k}$ is said to be {\it regular}
if it can be extended to a regular collection of $m+1$ complement subsets.

\begin{lemma}
Two distinct elements $U$ and $U'$ of ${\mathcal B}_{k}$ are adjacent
if and only if there exists a regular collection of complement subsets
${\mathcal R}^{1},\dots, {\mathcal R}^{m}\subset {\mathcal B}_{k}$ such that
$U$ and $U'$ belong to each ${\mathcal R}^{i}$, $i=1,\dots, m$.
\end{lemma}

\begin{proof}
It follows from the definition and Lemma 2.4.
\end{proof}

Now we can prove Theorem 2.5.
Suppose that $\Pi$ and $\Pi'$ satisfy the axiom ${\rm (P2)}$ and
consider an injection $f:{\mathcal G}_{k}(\Pi)\to{\mathcal G}_{k}(\Pi')$
sending base subsets to base subsets.

Let ${\mathcal B}_{k}$ be base subset of ${\mathcal G}_{k}(\Pi)$.
By our hypothesis, $f({\mathcal B}_{k})$ is a base subset of ${\mathcal G}_{k}(\Pi')$.
The mapping $f$ transfers inexact subsets of ${\mathcal B}_{k}$ to
inexact subsets of $f({\mathcal B}_{k})$
(indeed, if ${\mathcal R}$ is an inexact subset of ${\mathcal B}_{k}$
then there exists another base subset ${\mathcal B}_{k}'$
containing ${\mathcal R}$; since our mapping is injective,
the set $f({\mathcal R})$ is contained in the distinct
base subsets $f({\mathcal B}_{k})$ and $f({\mathcal B}_{k}')$,
this means that it is inexact).

We have established above that an inexact subset of
${\mathcal B}_{k}$ or $f({\mathcal B}_{k})$ is maximal if and only if
it consists of $\binom{n}{k+1}+\binom{n-1}{k-1}$ elements.
Hence maximal inexect subsets of ${\mathcal B}_{k}$
go over to maximal inexect subsets of $f({\mathcal B}_{k})$ and
$f$ sends complement subsets to complement subsets;
moreover, $f$ transfers regular collections of complement subsets
to regular collections (it follows from the definition).
Then Lemma 2.5 shows that the restriction of $f$ to ${\mathcal B}_{k}$
is adjacency preserving.
This implies the required, since for any two adjacent elements of
${\mathcal G}_{k}(\Pi)$ there exists a base subset containing them.

\subsection{Proof of Theorems 2.2 - 2.4.}
In this subsection we will always suppose that $n\ge 2k+1$ and the axiom (P2) holds
for both $\Pi$ and $\Pi'$.
Let $f$ be an injection of ${\mathcal G}_{k}(\Pi)$ to ${\mathcal G}_{k}(\Pi')$
sending base subsets to base subsets.
By Theorem 2.5, this mapping is adjacency preserving.
Then Proposition 2.2 shows that for any $(k-1)$-dimensional subspace $S$
of $\Pi$ there exists a $(k\pm 1)$-dimensional subspace $S'$ of $\Pi'$
such that
\begin{equation}
f({\mathcal G}_{k}(S))\subset {\mathcal G}_{k}(S').
\end{equation}
\begin{lemma}
For any $(k-1)$-dimensional subspace $S$ of $\Pi$
there is unique subspace $S'\in {\mathcal G}_{k\pm 1}(\Pi')$ satisfying the condition
$(2.5)$.
The dimension of this subspace may be equal to $k+1$ only for the case when $n=2k+1$.
\end{lemma}

\begin{proof}
Let ${\mathcal B}_{k-1}$ be a base subset  of ${\mathcal G}_{k-1}(\Pi)$ containing $S$
and $B$ be the base for $\Pi$ associated with ${\mathcal B}_{k-1}$.
Put ${\mathcal B}_{k}$ for the base subset of ${\mathcal G}_{k}(\Pi)$
defined by the base $B$;
this base subset will be called also associated with ${\mathcal B}_{k-1}$.
Then ${\mathcal B}'_{k}:=f({\mathcal B}_{k})$
is a base subset of ${\mathcal G}_{k}(\Pi')$.
For any number $m\in \{1,\dots, n-1\}$ we denote by ${\mathcal B}'_{m}$
the base subset of ${\mathcal G}_{m}(\Pi')$ associated with ${\mathcal B}'_{k}$.

It is trivial that
\begin{equation}
f({\mathcal G}_{k}(S)\cap {\mathcal B}_{k})\subset
{\mathcal G}_{k}(S')\cap {\mathcal B}'_{k}.
\end{equation}
Let us take two distinct subspaces $U_{1}$ and $U_{2}$ belonging to
${\mathcal G}_{k}(S)\cap {\mathcal B}_{k}$.
Then
$$S'=
\begin{cases}
\;f(U_{1})\cap f(U_{2})\;&\;\mbox{ if }\;\dim S'=k-1\\
\;\overline{f(U_{1})\cup f(U_{2})}\;&\;\mbox{ if }\;\dim S'=k+1.
\end{cases}
$$
The subspaces $f(U_{1})$ and $f(U_{2})$ belong to ${\mathcal B}'_{k}$,
hence $S'$ is an element of ${\mathcal B}'_{k\pm 1}$.
It is not difficult to see that ${\mathcal G}_{k}(S)\cap {\mathcal B}_{k}$
consists of $n-k+1$ elements and
$$|{\mathcal G}_{k}(S')\cap {\mathcal B}'_{k}|=
\begin{cases}
\;n-k+1\;&\;\mbox{ if }\;\dim S'=k-1\\
\;k+2\;&\;\mbox{ if }\; \dim S'=k+1\,.
\end{cases}
$$
The condition $2k+1 \le n$ implies the inequality
$$n-k+1\ge k+2$$
which can be replaced by the equality
only for the case when $n=2k+1$.
Since $f$ is injective, the inclusion (2.6) shows that
$$f({\mathcal G}_{k}(S)\cap {\mathcal B}_{k})=
{\mathcal G}_{k}(S')\cap {\mathcal B}'_{k}$$
and $S'$ may be $(k+1)$-dimensional only if $n=2k+1$.
Also this implies that $S'$ is uniquely defined.
\end{proof}

Denote by $f_{k-1}$ the mapping of ${\mathcal G}_{k-1}(\Pi)$
to $${\mathcal G}_{k-1}(\Pi')\cup {\mathcal G}_{k+1}(\Pi')$$
transferring each $(k-1)$-dimensional subspace $S$ of $\Pi$
to the subspace $S'\in {\mathcal G}_{k\pm 1}(\Pi')$
satisfying the condition (2.5).

\begin{lemma}
Suppose that $n=2k+1$ and there exists a subspace $S_{0}\in {\mathcal G}_{k-1}(\Pi)$
such that  $f_{k-1}(S_{0})$ is $(k+1)$-dimensional.
Then $f_{k-1}(S)$ is $(k+1)$-dimensional
for each $S\in {\mathcal G}_{k-1}(\Pi)$.
\end{lemma}

\begin{proof}
If $k=1$ then $S_{0}$ and $S$ are points; any two points are adjacent.
If $k>1$ then  $S_{0}$ and $S$ can be connected by a sequence of adjacent elements
of ${\mathcal G}_{k-1}(\Pi)$ (Proposition 2.1).
Thus we can restrict ourself only to the case when
$S_{0}$ and $S$ are adjacent.
For this case there exists a base subset ${\mathcal B}_{k-1}$ of ${\mathcal G}_{k-1}(\Pi)$
containing $S_{0}$ and $S$.
Let ${\mathcal B}_{k}$ be the base subset  of ${\mathcal G}_{k}(\Pi)$
associated with ${\mathcal B}_{k-1}$.
Define ${\mathcal B}'_{k}:=f({\mathcal B}_{k})$
and put ${\mathcal B}'_{m}$ for the base subset of ${\mathcal G}_{m}(\Pi')$
associated with ${\mathcal B}'_{k}$.

By the arguments given in the proof of Lemma 2.6, we have
$$f_{k-1}(S_{0})\in {\mathcal B}'_{k+1}\;\mbox{ and }\;
f_{k-1}(S)\in {\mathcal B}'_{k\pm 1};$$
besides
$$f({\mathcal G}_{k}(S_{0})\cap {\mathcal B}_{k})=
{\mathcal G}_{k}(f_{k-1}(S_{0}))\cap {\mathcal B}'_{k}$$
and
$$f({\mathcal G}_{k}(S)\cap {\mathcal B}_{k})=
{\mathcal G}_{k}(f_{k-1}(S))\cap {\mathcal B}'_{k}.$$
Since $S_{0}$ and $S$ are adjacent, the intersection of the sets
$${\mathcal G}_{k}(S_{0})\cap {\mathcal B}_{k}\; \mbox{ and }\;
{\mathcal G}_{k}(S)\cap {\mathcal B}_{k}$$
consists only of the $k$-dimensional subspace spanned by $S_{0}$ and $S$.
The mapping $f$ is injective and
\begin{equation}
{\mathcal G}_{k}(f_{k-1}(S_{0}))\cap{\mathcal G}_{k}(f_{k-1}(S))\cap {\mathcal B}'_{k}
\end{equation}
is a one-element set.
If the dimension of $f_{k-1}(S)$ is equal to $k-1$ then
the set (2.7) is not empty only for the case when
$$f_{k-1}(S)\subset f_{k-1}(S_{0}).$$
However, if this inclusion holds then (2.7) contains
more than one element.
This means that $f_{k-1}(S)$ is $(k+1)$-dimensional.
\end{proof}

Thus $f_{k-1}$ is a mapping of ${\mathcal G}_{k-1}(\Pi)$
to ${\mathcal G}_{k-1}(\Pi')$ or ${\mathcal G}_{k+1}(\Pi')$
such that the equality
$$ f({\mathcal G}_{k}(S))\subset{\mathcal G}_{k}(f_{k-1}(S))$$
holds for any $(k-1)$-dimensional subspace $S$ of $\Pi$.
If $f_{k-1}$ is a mapping to ${\mathcal G}_{k-1}(\Pi')$ then
\begin{equation}
f_{k-1}({\mathcal G}_{k-1}(U))\subset {\mathcal G}_{k-1}(f(U))
\;\;\;\;\;\;\forall\;U\in{\mathcal G}_{k}(\Pi).
\end{equation}
Note that $f_{k-1}$ may be a mapping to ${\mathcal G}_{k+1}(\Pi')$ only for the case when
$n=2k+1$.

\begin{lemma}
The mapping $f_{k-1}$ transfers base subsets to base subsets.
\end{lemma}

\begin{proof}
Let ${\mathcal B}_{k-1}$ be a base subset of ${\mathcal G}_{k-1}(\Pi)$
and ${\mathcal B}_{k}$ be the associated base subset of ${\mathcal G}_{k}(\Pi)$.
Then $f_{k-1}({\mathcal B}_{k-1})$ is the base subset of
${\mathcal G}_{k\pm 1}(\Pi')$ associated with $f({\mathcal B}_{k})$.
\end{proof}

\begin{lemma}
Suppose that $f_{k-1}$ is induced by some strong embedding $g$ of $\Pi$ to $\Pi'$.
Then $f=G_{k}(g)$.
\end{lemma}

\begin{proof}
Since $f_{k-1}$ is induced by $g$,
for any $k$-dimensional subspace $U$ of $\Pi$ we have
$$f_{k-1}({\mathcal G}_{k-1}(U))\subset {\mathcal G}_{k-1}(\overline{g(U)})$$
and $f_{k-1}$ is injective.
The last equality together with (2.8) show that
$$f_{k-1}({\mathcal G}_{k-1}(U))\subset
{\mathcal G}_{k-1}(f(U))\cap {\mathcal G}_{k-1}(\overline{g(U)}).$$
If $f(U)$ and $\overline{g(U)}$ are different then the set
$${\mathcal G}_{k-1}(f(U))\cap {\mathcal G}_{k-1}(\overline{g(U)})$$
contains at most one element.
This contradicts to the injectivity of $f_{k-1}$.
Therefore, $f(U)=\overline{g(U)}$ and we get the required.
\end{proof}

\begin{proof}[Proof of Theorem 2.3]
Let $f$ be an injection of ${\mathcal L}$ to ${\mathcal L}'$
sending base subsets to base subsets and $n>3$.
Then $f$ induces certain mapping $f_{0}:P\to P'$.
Lemma 2.8 says that $f_{0}$ maps bases for $\Pi$ to bases for $\Pi'$;
in particular, $f_{0}$ is injective (see Subsection 1.4, the proof of Theorem 1.3).
Moreover, (2.8) guarantees that $f_{0}$ is collinearity preserving.
Therefore, $f_{0}$ is a strong embedding. By Lemma 2.9,
this strong embedding induces $f$.
\end{proof}

\begin{proof}[Proof of Theorem 2.2]
Let $\Pi$ and $\Pi'$ be projective and
$f$ be a mapping of ${\mathcal G}_{k}(\Pi)$ to ${\mathcal G}_{k}(\Pi')$
which transfers base subsets to base subsets (we do not require the injectivity).
It was noted above (Subsection 2.6) that we need to prove Theorem 2.2 only
for the case when $n\ge 2k+1$.

Since for any two elements of ${\mathcal G}_{k}(\Pi)$
there exists a base subset containing them, the mapping $f$ is injective
(suppose that $S$ and $U$ are distinct elements of ${\mathcal G}_{k}(\Pi)$
such that $f(S)=f(U)$ and ${\mathcal B}_{k}$ is a base subset
containing $S$ and $U$, then the set $f({\mathcal B}_{k})$ has less than
$\binom{n+1}{k+1}$ elements, this means that $f({\mathcal B}_{k})$ is not a base subset).
Therefore, $f$ is adjacency preserving (Theorem 2.5)
and induces certain  mapping $f_{k-1}$ of
${\mathcal G}_{k-1}(\Pi)$ to ${\mathcal G}_{k-1}(\Pi')$ or ${\mathcal G}_{k+1}(\Pi')$.

Now we prove by induction that for the fist case
($f_{k-1}$ is a mapping to ${\mathcal G}_{k-1}(\Pi')$)
$f$ is induced by a strong embedding of $\Pi$ to $\Pi'$.
It was established only for $k=1$ (Theorem 2.3). Suppose that $k\ge 2$.
The mapping $f_{k-1}$ transfers base subsets to base subsets (Lemma 2.8)
and the arguments given above show that $f_{k-1}$ is injective.
We have
$$2(k-1)+1<2k+1\le n$$
and $f_{k-1}$ induces a mapping of ${\mathcal G}_{k-2}(\Pi)$ to ${\mathcal G}_{k-2}(\Pi')$.
The inductive hypothesis implies the existence of
a strong embedding of $g$ of $\Pi$ to $\Pi'$ such that $G_{k-1}(g)=f_{k-1}$.
By Lemma 2.9, the embedding $g$ induces $f$.

If $f_{k-1}$ is a mapping to ${\mathcal G}_{k+1}(\Pi')$
then $n=2k+1$.
For this case $f$ can be considered as an
injection of ${\mathcal G}_{k}(\Pi)$ to ${\mathcal G}_{k}(\Pi'^{*})$
sending bases subsets to bases subsets and inducing a mapping
of ${\mathcal G}_{k-1}(\Pi)$ to ${\mathcal G}_{k-1}(\Pi'^{*})$.
This implies that $f$ is induced by a strong embedding of $\Pi$ to $\Pi'^{*}$.
\end{proof}

\begin{proof}[Proof of Theorem 2.4]
Suppose that $f:{\mathcal G}_{k}(\Pi)\to {\mathcal G}_{k}(\Pi')$
is a bijection such that $f$ and $f^{-1}$ map base subsets to base subsets
and $n>2k+1$.
Then $f_{k-1}$ is a mapping of
${\mathcal G}_{k-1}(\Pi)$ to ${\mathcal G}_{k-1}(\Pi')$
and $f^{-1}$ induces the mapping of ${\mathcal G}_{k-1}(\Pi')$
to ${\mathcal G}_{k-1}(\Pi)$ which is inverse to $f_{k-1}$.
Hence $f_{k-1}$ is bijective.
By Lemma 2.8, $f_{k-1}$ transfers base subsets to base subsets
and the same holds for the inverse mapping.
So, we can prove by induction that $f$ is induced by
a collineation of $\Pi$ to $\Pi'$.
\end{proof}

\end{document}